# A Robust Consensus Algorithm for Current Sharing and Voltage Regulation in DC Microgrids

Michele Cucuzzella, Sebastian Trip, Claudio De Persis, Xiaodong Cheng, Antonella Ferrara, and Arjan van der Schaft

*Abstract*—In this paper a novel distributed control algorithm for current sharing and voltage regulation in Direct Current (DC) microgrids is proposed. The DC microgrid is composed of several Distributed Generation units (DGUs), including Buck converters and current loads. The considered model permits an arbitrary network topology and is affected by unknown load demand and modelling uncertainties. The proposed control strategy exploits a communication network to achieve proportional current sharing using a consensus-like algorithm. Voltage regulation is achieved by constraining the system to a suitable manifold. Two robust control strategies of Sliding Mode (SM) type are developed to reach the desired manifold in a finite time. The proposed control scheme is formally analyzed, proving the achievement of proportional current sharing, while guaranteeing that the weighted average voltage of the microgrid is identical to the weighted average of the voltage references.

*Index Terms*—DC Microgrids, Sliding mode control, Uncertain systems, Current sharing, Voltage regulation.

## I. INTRODUCTION

IN the last decades, due to economic, technological and environmental aspects, the main trends in power systems focused on the modification of the traditional power generation and transmission systems towards incorporating smaller Distributed Generation units (DGUs). Moreover, the ever-increasing energy demand and the concern about the climate change have encouraged the wide diffusion of Renewable Energy Sources (RES). The so-called microgrids have been proposed as conceptual solutions to integrate different types of RES and to electrify remote areas. Microgrids are low-voltage electrical distribution networks, composed of clusters of DGUs, loads and storage systems interconnected through power lines [2].

Due to the widespread use of Alternate Current (AC) electricity in most industrial, commercial and residential applications, the recent literature on this topic mainly focused on

This work is supported by the EU Project 'MatchIT' (project number: 82203). Also, this work is part of the research programme ENBARK+ with project number 408.urs+.16.005, which is (partly) financed by the Netherlands Organisation for Scientific Research (NWO). Preliminary results have appeared in [1]

M. Cucuzzella (corresponding author), S. Trip, C. De Persis and X. Cheng are with Jan C. Wilems Center for Systems and Control, ENTEG, Faculty of Science and Engineering, University of Groningen, Nijenborgh 4, 9747 AG Groningen, the Netherlands, (email: {m.cucuzzella, s.trip, c.de.persis, x.cheng}@rug.nl).

A. Ferrara is with the Dipartimento di Ingegneria Industriale e dell'Informazione, University of Pavia, via Ferrata 5, 27100 Pavia, Italy, (e-mail: antonella.ferrara@unipv.it).

A. van der Schaft is with the Johann Bernoulli Institute for Mathematics and Computer Science, University of Groningen, Nijenborgh 9, 9747 AG Groningen, the Netherlands, (email: a.j.van.der.schaft@rug.nl).

AC microgrids [3]–[7]. However, several sources and loads (e.g. photovoltaic panels, batteries, electronic appliances and electric vehicles) can be directly connected to DC microgrids by using DC-DC converters. Indeed, several aspects make DC microgrids more efficient and reliable than AC microgrids [8]: *i*) lossy DC-AC and AC-DC conversion stages are reduced, *ii*) there is not reactive power, *iii*) harmonics are not present, *iv*) frequency synchronization is overcame, *v*) the skin effect is absent. Moreover, a DC microgrid can be connected to an islanded AC microgrid (even to the main grid) by a DC-AC bidirectional converter, forming a so-called hybrid microgrid [9]. Moreover, the growing need of interconnecting distant power networks (e.g. off-shore wind farms) has encouraged the use of High Voltage Direct Current (HVDC) technology [10]–[12], which is advantageous not only for long distances, but also for underwater cables, asynchronous networks and grids running at different frequencies [13]. Finally, DC microgrids are widely deployed in aircrafts and trains, and recently used in modern design for ships and large charging facilities for electric vehicles. For all these reasons, DC microgrids are attracting growing interest and receive much research attention.

### A. Literature review

Two main control objectives in DC microgrids are voltage regulation and current sharing (or, equivalently, load sharing). Regulating the voltages is required to ensure a proper functioning of connected loads [14]–[17], whereas current sharing prevents the overstressing of any source. Moreover, since a microgrid can include DGUs with different generation capacity, it is often desired in practical cases that the DGUs share the total current demand proportionally to their generation capacity. In order to achieve both objectives, hierarchical control schemes are conventionally adopted [18]. Generally, the requirement of current sharing does not permit to regulate the voltage at each node towards the corresponding desired value. Then, a reasonable alternative is to satisfy the voltage requirement defined in [19], according to which the average voltage across the whole microgrid (not a specific node) should be regulated at the global voltage set point (e.g., the average of the voltage references). This kind of voltage regulation is called *global voltage regulation* or *voltage balancing* (see for instance [20]–[24] and the references therein).

In the literature, these control problems in DC microgrids have been addressed by different control approaches and schemes, and we discuss a few of them. To compensate the voltage steady state error due to primary droop controller, a distributed secondary controller based on averaging the total current supplied by the sources is proposed in [25]. Yet, for the stability analysis, fast dynamics are neglected and only the small-signal model is considered. Distributed secondary integral control strategies that are able to achieve proportional load sharing and voltage regulation are formally analyzed in [26], neglecting inductive lines. In [19] each power converter is equipped with current and voltage regulators in order to achieve both proportional load sharing and voltage regulation. However, the achievement of voltage regulation requires the use of an observer to estimate the global average voltage, leading to more complicated controller implementations. In [24] the authors propose a consensus-based secondary controller for proportional current sharing and global voltage regulation for resistive networks. However, proportional current sharing is achieved under the restrictive assumptions that the line resistances are known and the electrical and communication graphs are identical. A consensus algorithm that guarantees power sharing in presence of 'ZIP' (constant impedance, constant current, constant power) loads, as well as preservation of the weighted geometric average of the source voltages is designed and formally analyzed in [27]. However, only pure resistive networks are considered and the steady state voltages strongly depend on the voltage initial conditions.

### B. Main contributions

This paper proposes a novel robust control algorithm to obtain simultaneously proportional current sharing among the DGUs and a form of voltage regulation in the DC power network, where the interconnecting lines of the microgrid are assumed to be resistive-inductive. In order to achieve current sharing, a communication network is exploited where each DGU communicates in real-time the value of its generated current to its neighbouring DGUs. Adding this additional communication layer to achieve current sharing, leading to a distributed controller, has been widely adopted and studied thoroughly. In comparison to the existing results in the literature, we additionally propose the design of a manifold that couples the aforementioned objective of current sharing to the objective of voltage regulation. By doing this, the proposed control algorithm guarantees that the weighted average voltage of the microgrid is equal to the weighted average of the reference voltages, where the weights depend on the DGUs generation capacities, performing the so called global voltage regulation or voltage balancing [19], [24]. This is achieved independently of the initial voltage conditions, facilitating Plug-and-Play capabilities.

To constrain the state of the system to the designed manifold in a finite time, we propose robust controllers of Sliding Mode (SM) type [28], [29]. SM control is appreciated for its robustness property against a wide class of modelling uncertainties and external disturbances, commonly present in DC microgrids. In this paper, we first propose a Second Order Sliding Mode (SOSM) controller that determines the, possibly non-constant, switching frequency of the power converter, which might lead increased the power losses. Then, to overcome this issue, we additionally propose a third order sliding mode controller (3SM) to obtain a continuous control signal that can be used as the duty cycle of the power converter. Furthermore, the proposed control solution is robust with respect to failed communication. In fact, if the communication among the DGUs is disabled, then the voltage of each node converges in a finite time to the corresponding reference value. For the considered microgrid model, convergence to the state of current sharing and voltage regulation is theoretically analyzed, and we show that convergence is achieved globally, for any initialization of the microgrid.

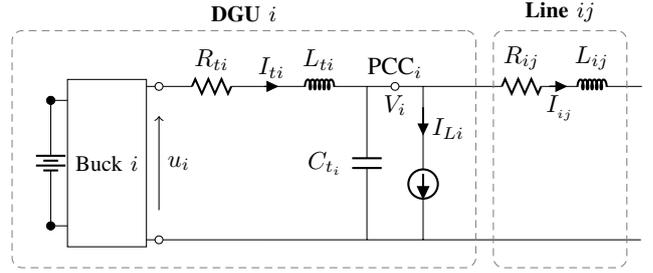

Fig. 1. Electrical scheme of DGU $i$ and line $ij$.

### C. Outline

The remainder of this paper is organized as follows. In Section II the microgrid model is presented, while in Section III the control problem is formulated. In Section IV the proposed manifold-based consensus algorithm is designed, and in Section V sliding mode control strategies are proposed to reach the desired manifold. In Section VI the stability properties of the controlled system are analyzed, while in Section VII the simulation results are illustrated and discussed. Some conclusions are gathered in Section VIII.

## II. DC MICROGRID MODEL

In this work we consider a typical Buck converter-based DC microgrid of which a schematic electrical diagram is provided in Figure 1. By applying the Kirchhoff's current (KCL) and voltage (KVL) laws, the governing dynamic equations of the $i$-th node (DGU) are the following:

$$
\begin{aligned}
L_{ti}\dot{I}_{ti} &= -R_{ti}I_{ti} - V_i + u_i \\
C_{ti}\dot{V}_i &= I_{ti} - I_{Li} - \sum_{j \in \mathcal{N}_i} I_{ij},
\end{aligned}
\tag{1}
$$

where $\mathcal{N}_i$ is the set of nodes (i.e., the DGUs) connected to the $i$-th DGU by distribution lines, while the control input $u_i$ represents the Buck converter output voltage*. The current from DGU $i$ to DGU $j$ is denoted by $I_{ij}$, and its dynamic is given by

$$
L_{ij}\dot{I}_{ij} = (V_i - V_j) - R_{ij}I_{ij}.
\tag{2}
$$

The symbols used in (1) and (2) are described in Table I.

---

*Note that $u_i$ in (1) can be expressed as $\delta_i V_{DC_i}$, where $\delta_i$ is the duty cycle of the Buck converter $i$ and $V_{DC_i}$ is the DC voltage source provided by a generic renewable energy source or a battery at node $i$.

| State variables | |
|---|---|
| $I_{ti}$ | Generated current |
| $V_i$ | Load voltage |
| $I_{ij}$ | Exchanged current |
| **Parameters** | |
| $R_{ti}$ | Filter resistance |
| $L_{ti}$ | Filter inductance |
| $C_{ti}$ | Shunt capacitor |
| $R_{ij}$ | Line resistance |
| $L_{ij}$ | Line inductance |
| **Inputs** | |
| $u_i$ | Control input |
| $\bar{I}_{ti}$ | Unknown current demand |

The overall network is represented by a connected and undirected graph $\mathcal{G} = (\mathcal{V}, \mathcal{E})$, where the nodes, $\mathcal{V} = \{1, ..., n\}$, represent the DGUs and the edges, $\mathcal{E} = \{1, ..., m\}$, represent the distribution lines interconnecting the DGUs. The network topology is represented by its corresponding incidence matrix $\mathcal{B} \in \mathbb{R}^{n \times m}$. The ends of edge $k$ are arbitrarily labeled with a $+$ and a $-$, and the entries of $\mathcal{B}$ are given by

$$\mathcal{B}_{ik} = \begin{cases} +1 & \text{if } i \text{ is the positive end of } k \\ -1 & \text{if } i \text{ is the negative end of } k \\ 0 & \text{otherwise.} \end{cases}$$

Consequently, the overall microgrid system can be written compactly for all nodes $i \in \mathcal{V}$ as

$$\begin{aligned} L_t \dot{I}_t &= -R_t I_t - V + u \\ C_t \dot{V} &= I_t + \mathcal{B}I - I_L \\ L \dot{I} &= -\mathcal{B}^T V - RI, \end{aligned} \quad (3)$$

where $I_t, V, I_L, u \in \mathbb{R}^n$, and $I \in \mathbb{R}^m$. Moreover, $C_t, L_t, R_t \in \mathbb{R}^{n \times n}$ and $R, L \in \mathbb{R}^{m \times m}$ are positive definite diagonal matrices, e.g. $R_t = \text{diag}(R_{t_1}, \ldots, R_{t_n})$. To permit the controller design in the next sections, the following assumption is introduced on the available information of the system:

*Assumption 1:* (**Available information**) The state variables $I_{ti}$ and $V_i$ are locally available at the $i$-th DGU. The network parameters $R_t, L_t, C_t, R, L$ and the current demand $I_L$ are constant and unknown, but with known bounds.

*Remark 1:* (**Varying parameters and current demand**) We assume that the parameters and the current demand are constant, to allow for a steady state solution and to theoretically analyze the stability of the microgrid. Yet, the control strategy that we propose in the next sections is applicable even if this assumption is removed.

*Remark 2:* (**Kron reduction**) Note that in (1), the load currents are located at the PCC of each DGU (see also Figure 1). This situation is generally obtained by a Kron reduction of the original network, yielding an equivalent representation of the network [26]. It is important to realize that the network (topology) of the Kron reduced network is generally unknown and differs from the original network. It is therefore desirable that a control structure is independent of the underlying distribution network.

## III. CURRENT SHARING AND VOLTAGE BALANCING

In this section we make the considered control objectives explicit. First, we note that for a given constant control input $\bar{u}$, a steady state solution $(\bar{I}_t, \bar{V}, \bar{I})$ to system (3) satisfies

$$\begin{aligned} \bar{V} &= -R_t \bar{I}_t + \bar{u} \\ -\mathcal{B}\bar{I} &= \bar{I}_t - I_L \\ \bar{I} &= -R^{-1}\mathcal{B}^T \bar{V}. \end{aligned} \quad (4)$$

The second line of (4) implies[†] that at steady state the total generated current $\mathbb{1}_n^T \bar{I}_t$ is equal to the total current demand $\mathbb{1}_n^T I_L$. To improve the generation efficiency, it is generally desired that the total current demand is shared among the various DGUs proportionally to the generation capacity of their corresponding energy sources (proportional current sharing). This desire can be expressed as $w_i \bar{I}_{ti} = w_j \bar{I}_{tj}$ for all $i, j \in \mathcal{V}$, where $w_i$ relates to the generation capacity of converter $i$, and leads to the first objective concerning the desired steady state value of the generated currents $\bar{I}_t$.

*Objective 1:* (**Proportional Current sharing**)

$$\lim_{t \to \infty} I_t(t) = \bar{I}_t = W^{-1} \mathbb{1}_n i_t^\star, \quad (5)$$

with $i_t^\star = \mathbb{1}_n^T I_L / (\mathbb{1}_n^T W^{-1} \mathbb{1}_n) \in \mathbb{R}$, $W = \text{diag}\{w_1, \ldots, w_n\}$, $w_i > 0$, for all $i \in \mathcal{V}$.

Note that (5) indeed satisfies $\mathbb{1}_n^T \bar{I}_t = \mathbb{1}^T W^{-1} \mathbb{1}_n i_t^\star = \mathbb{1}_n^T I_L$. From the second and third lines of (4) it follows that the corresponding steady state voltages $\bar{V}$ satisfy $\mathcal{B}R^{-1}\mathcal{B}^T \bar{V} = W^{-1} \mathbb{1}_n i_t^\star - I_L$, that prescribes the value of the required differences in voltages, $\mathcal{B}^T \bar{V}$, achieving proportional current sharing. This admits the freedom to shift all steady state voltages with the same constant value, since $\mathcal{B}^T \bar{V} = \mathcal{B}^T (\bar{V} + a\mathbb{1}_n)$, with $a \in \mathbb{R}$ any scalar. To define the optimal steady state voltages, we assume that for every DGU $i$, there exists a desired reference voltage $V_i^\star$.

*Assumption 2:* (**Desired voltages**) There exists a constant reference voltage $V_i^\star$ at the PCC, for all $i \in \mathcal{V}$.

Often the values for $V_i^\star$ are chosen identical for all $i \in \mathcal{V}$, and are set to the desired voltage level of the overall network. Generally, the requirement of current sharing does not permit for $\bar{V} = V^\star$, and might cause voltages deviations from the corresponding reference values. Then, a reasonable alternative is to keep the weighted average value of the PCC voltages at the steady state identical to the weighted average value of the desired reference voltages $V^\star$ (voltage balancing) [24]. Particularly, we choose the weights to be $1/w_i$, for all $i \in \mathcal{V}$, such that at the converters with a relatively large generation capacity, there is a relatively small voltage deviation. It is indeed a standard practise that the sources with the largest generation capacity determine the grid voltage. Therefore, given a $V^\star$, we aim at designing a controller that, in addition to Objective 1, also guarantees voltage balancing, i.e.,

---

[†]The incidence matrix $\mathcal{B}$, satisfies $\mathbb{1}_n^T \mathcal{B} = \mathbf{0}$, where $\mathbb{1}_n \in \mathbb{R}^n$ is the vector consisting of all ones.

*Objective 2:* (**Voltage balancing**)

$$\lim_{t \to \infty} \mathbb{1}_n^T W^{-1} V(t) = \mathbb{1}_n^T W^{-1} \overline{V} = \mathbb{1}_n^T W^{-1} V^\star. \quad (6)$$

*Remark 3:* (**Equal current sharing**) Note that by setting in (5) and (6) the weights $w_i$, for all $i \in \mathcal{V}$, identical, the total current demand is equally shared among the DGUs and the arithmetic average of the microgrid voltage is equal to the arithmetic average of the voltage references.

By substituting (5) and (6), in (4), one can easily verify that achieving Objective 1 and Objective 2 prescribes the (optimal) steady state output voltages of the Buck converters, $u = \overline{u}^{opt}$.

*Lemma 1:* (**Optimal feedforward input**) If system (3), at steady state, achieves Objective 1 and Objective 2, then the control input $u$ to system (3) is given by

$$\overline{u}^{opt} = -\left(\mathcal{B}R^{-1}\mathcal{B}^T - \Psi\right)^{-1}\left(\Psi V^\star + I_L\right), \quad (7)$$

with

$$\Psi = \frac{(\mathbb{1}_n + \mathcal{B}R^{-1}\mathcal{B}^T R_t)W^{-1}\mathbb{1}_n\mathbb{1}_n^T W^{-1}}{\mathbb{1}_n^T W^{-1} R_t W^{-1}\mathbb{1}_n}, \quad (8)$$

and $\mathbb{1}_n \in \mathbb{R}^{n \times n}$ the identity matrix.

*Proof:* When Objective 1 and Objective 2 hold, the steady state of (3) necessarily satisfies

$$\begin{aligned}
\mathbf{0} &= -R_t W^{-1}\mathbb{1}_n i_t^\star - \overline{V} + \overline{u}^{opt} \\
\mathbf{0} &= W^{-1}\mathbb{1}_n i_t^\star - \mathcal{B}R^{-1}\mathcal{B}^T \overline{V} - I_L \\
\mathbf{0} &= \mathbb{1}_n^T W^{-1} \overline{V} - \mathbb{1}_n^T W^{-1} V^\star,
\end{aligned} \quad (9)$$

with $i_t^\star = \mathbb{1}_n^T I_L / (\mathbb{1}_n^T W^{-1}\mathbb{1}_n) \in \mathbb{R}$. A tedious, but straightforward, calculation permits to solve (9) for $\overline{u}^{opt}$, yielding (7). ∎

In order to determine (7), exact knowledge of almost all network parameters, as well as the current demand $I_L$, is required. Since this information is not available (see also Assumption 1), we propose in the next sections distributed controllers that, provably, achieve voltage balancing using only local measurements of $V_i$, and that achieve proportional current sharing by exchanging information on $I_{ti}$ among neighbours over a communication network. In the remainder of this section we further elaborate on the steady state voltages imposed by the control objectives.

### A. Steady state voltages

First, we notice that it follows from (5) and (9) that the steady state voltages $\overline{V}$ satisfy

$$\overline{V} = -\frac{R_t W^{-1}\mathbb{1}_n\mathbb{1}_n^T I_L}{\mathbb{1}_n^T W^{-1}\mathbb{1}_n} + \overline{u}^{opt}. \quad (10)$$

From (7) and (10) it is evident that the steady state values of the voltages at each node depend on the loads $I_L$ and the voltage references $V^\star$. Since $V^\star$ is free to design, it can be potentially chosen in such a way that too low or too high voltages are avoided. To help the design of $V^\star$, we show that the the steady state voltages $\overline{V}_i$, for all $i \in \mathcal{V}$, are shifted by the same quantity, when $V^\star$ is altered.

*Lemma 2:* (**Voltage shifting property**) Let Objective 1 and Objective 2 hold, and let $\overline{V}_{(1)} \in \mathbb{R}^n$ denote the steady state voltage value associated to the voltage reference $V_{(1)}^\star \in \mathbb{R}^n$. Consider the new voltage reference $V_{(2)}^\star \in \mathbb{R}^n$ and the corresponding steady state voltage value $\overline{V}_{(2)} \in \mathbb{R}^n$. Then, $\Delta \overline{V} = \overline{V}_{(2)} - \overline{V}_{(1)}$ satisfies

$$\Delta \overline{V} = \mathbb{1}_n \frac{\mathbb{1}_n^T W^{-1}\Delta V^\star}{\mathbb{1}_n^T W^{-1}\mathbb{1}_n}, \quad (11)$$

with $\Delta V^\star = V_{(2)}^\star - V_{(1)}^\star$.

*Proof:* When Objective 2 holds, we have

$$\mathbb{1}_n^T W^{-1}(\overline{V}_{(1)} + \Delta \overline{V}) = \mathbb{1}_n^T W^{-1}(V_{(1)}^\star + \Delta V^\star), \quad (12)$$

which implies $\mathbb{1}_n^T W^{-1}\Delta \overline{V} = \mathbb{1}_n^T W^{-1}\Delta V^\star$. Bearing in mind that the voltage differences between any node of the microgrid are prescribed by the achievement of current sharing (see the paragraph below Objective 1), we have $\mathcal{B}^T \overline{V}_{(1)} = \mathcal{B}^T \overline{V}_{(2)}$, implying $\Delta \overline{V} = \mathbb{1}_n \nu$, with $\nu = \mathbb{1}_n^T W^{-1}\Delta V^\star / \mathbb{1}_n^T W^{-1}\mathbb{1}_n$, i.e., all the voltages are shifted by the same quantity. ∎

Consequently, any node $i$ in the network can lower or increase its steady state voltage $\overline{V}_i$, by adjusting its own reference $V_i^\star$. Although, the design and the analysis of a voltage reference generator is postponed to a future research, the property proven in Lemma 2 could be exploited to tune the references in order to avoid that the voltages at some nodes are lower or higher than some given thresholds. ∎

## IV. A Manifold-Based Consensus Algorithm

In this section we introduce the key aspects of the proposed solution to simultaneously achieve Objective 1 and Objective 2, consisting of a consensus algorithm and the design of a manifold to where the solutions to the system should converge. First, we augment system (3) with additional state variables (distributed integrators) $\theta_i$, $i \in \mathcal{V}$, with dynamics given by

$$\dot{\theta}_i = -\sum_{j \in \mathcal{N}_i^c} \gamma_{ij}(w_i I_{ti} - w_j I_{tj}), \quad (13)$$

where $\mathcal{N}_i^c$ is the set of the DGUs that communicate with the $i$-th DGU, $\gamma_{ij} = \gamma_{ji} \in \mathbb{R}_{>0}$ are additional gain constants, and $w_i, w_j \in \mathbb{R}_{>0}$ are constant weights depending on the DGUs generation capacity. Let $\mathcal{L}_c$ denote the (weighted) Laplacian matrix associated with the communication graph, which can be different from the topology of the (reduced) microgrid. Then, the dynamics in (13) can be expressed compactly for all nodes $i \in \mathcal{V}$ as

$$\dot{\theta} = -\mathcal{L}_c W I_t, \quad (14)$$

that indeed has the form of a consensus protocol, permitting a steady state where $W\overline{I}_t \in \text{im}(\mathbb{1}_n)$ (see also Objective 1). We impose the following restrictions on (14):

*Assumption 3:* (**Controller structure**) For all $i \in \mathcal{V}$, the integrators states $\theta_i$ are initialized such that $\mathbb{1}_n^T \theta(0) = 0$. Furthermore, the graph corresponding to the topology of the communication network is undirected and connected. ∎

The most straightforward choice of initialization of the state $\theta_i(0)$, that satisfies Assumption 3, is to initialize all $\theta_i$ to zero, i.e. $\theta(0) = \mathbf{0}$. Whereas connectedness of the communication graph is needed to ensure current sharing among *all* DGUs, the consequence of the required initialization of $\theta$ is that the average value of the entries of $\theta$ is preserved and identical to zero for all $t \geq 0$, as proved in the following lemma:

*Lemma 3:* (**Preservation of** $\mathbb{1}_n^T \theta$) Let Assumption 3 hold. Given system (14), the average value $\frac{1}{n} \sum_{i \in \mathcal{V}} \theta_i$ is preserved, i.e.,

$$\frac{1}{n} \mathbb{1}_n^T \theta(t) = \frac{1}{n} \mathbb{1}_n^T \theta(0) \quad \text{for all } t \geq 0. \tag{15}$$

*Proof:* Pre-multiplying both sides of (14) by $\mathbb{1}_n^T$ yields

$$\mathbb{1}_n^T \dot{\theta} = -\mathbb{1}_n^T \mathcal{L}_c W I_t = 0, \tag{16}$$

where $\mathbb{1}_n^T \mathcal{L}_c = \mathbf{0}$, follows from $\mathcal{L}_c$ being the Laplacian matrix associated with an undirected graph. ∎

The fact that $\mathbb{1}_n^T \theta(t) = 0$, is essential to the second aspect of the proposed solution, the design of a manifold. Bearing in mind Objective 2, we propose the following desired manifold:

$$\{(I_t, V, I, \theta) : W^{-1}(V - V^\star) - \theta = \mathbf{0}\}. \tag{17}$$

Indeed, exploiting the preservation of $\mathbb{1}_n^T \theta$, we have on the desired manifold (17), $\mathbb{1}_n^T W^{-1} V = \mathbb{1}_n^T (\theta + W^{-1} V^\star) = \mathbb{1}_n^T W^{-1} V^\star$. Constraining the solutions to a system to a specific manifold is typical for sliding mode based controllers, and we will discuss some suitable controller designs in the next section.

*Remark 4:* (**Plug-and-Play**) The main results in this work assume a constant network topology. Nevertheless, an interesting extension is to consider the plugging in or out of various converters. The analysis of the corresponding switched/hybrid system is outside the scope of this work. Here, we merely describe how the required initialization $\theta_i$ should be extended towards the setting of changing topologies, in order to preserve the crucial property $\mathbb{1}_n^T \theta = 0$. First, if a new DGU (say DGU$_{n+1}$) wants to join the network, its integrator state is initialized to zero, i.e., $\theta_{n+1}(t_{new}) = 0$, $t_{new}$ being the time instant when DGU$_{n+1}$ is plugged-in. Second, if a DGU (say DGU $i$) is unplugged at the time instant $t_{out}$, we let $\theta_i(t) = \theta_i(t_{out})$ for all $t > t_{out}$, without re-setting any integrator. If DGU $i$ wants to join again the network at the time instant $t_{in} > t_{out}$, the dynamic of $\theta_i$ is described again by (13) for all $t > t_{in}$. Since $\theta_i(t_{in}) = \theta_i(t_{out})$, also the plug-in operation occurs without re-setting any integrator state.

## V. SLIDING MODE CONTROLLERS

We now propose a Distributed Second Order Sliding Mode (D-SOSM) control law, and a Distributed Third Order Sliding Mode (D-3SM) control law, to steer, in a finite time, the state of system (3), augmented with (14), to the desired manifold (17). As will be discussed in the coming subsections, the choice of the particular control law, D-SOSM or D-3SM, depends on the desired implementation.

First, to facilitate the upcoming discussion, we recall the following definitions that are essential to sliding mode control:

*Definition 1:* (**Sliding function**) Consider system

$$\dot{x} = \zeta(x, u), \tag{18}$$

with state $x \in \mathbb{R}^n$, and input $u \in \mathbb{R}^m$. The *sliding function* $\sigma(x) : \mathbb{R}^n \to \mathbb{R}^m$ is a sufficiently smooth output function of system (18).

*Definition 2:* ($r$**−sliding manifold**) The $r$−*sliding manifold*[‡] is given by

$$\{x \in \mathbb{R}^n, u \in \mathbb{R}^m : \sigma = L_\zeta \sigma = \cdots = L_\zeta^{(r-1)} \sigma = \mathbf{0}\}, \tag{19}$$

where $L_\zeta^{(r-1)} \sigma(x)$ is the $(r-1)$-th order Lie derivative of $\sigma(x)$ along the vector field $\zeta(x, u)$. With a slight abuse of notation we also write $L_\zeta \sigma(x) = \dot{\sigma}(x)$, and $L_\zeta^{(2)} \sigma(x) = \ddot{\sigma}(x)$.

*Definition 3:* ($r$**−order sliding mode (controller)**) A $r$−*order sliding mode* is enforced from $t = T_r \geq 0$, when, starting from an initial condition, the state of (18) reaches the $r$−sliding manifold, and remains there for all $t \geq T_r$. The order of a *sliding mode controller* is identical to the order of the sliding mode that it is aimed at enforcing.

Bearing in mind the definitions above and the desired manifold (17), we consider the following sliding function $\sigma \in \mathbb{R}^n$:

$$\sigma(V, \theta) = W^{-1}(V - V^\star) - \theta. \tag{20}$$

### A. Second order SM control: variable switching frequency

Regarding the sliding function (20) as the output function of system (3), (14), it appears that the relative degree[§] is two. This implies that a second order sliding mode (SOSM) controller can be *naturally* applied in order to make the state of the controlled system reach, in a finite time, the sliding manifold $\{(I_t, V, I, \theta) : \sigma = \dot{\sigma} = \mathbf{0}\}$. According to the SOSM control theory, the auxiliary variables $\xi_1 = \sigma$ and $\xi_2 = \dot{\sigma}$ have to be defined, resulting in the so-called auxiliary system

$$\begin{aligned} \dot{\xi}_1 &= \xi_2 \\ \dot{\xi}_2 &= b(I_t, V, I, u) + G_d u. \end{aligned} \tag{21}$$

Taking into account the expressions for $\sigma$ and $\dot{\sigma}$, a straightforward calculation shows that, in the auxiliary system (21), the expression for $b \in \mathbb{R}^n$ is given by

$$\begin{aligned} b = &-\left(W^{-1} C_t^{-1} + \mathcal{L}_c W\right) L_t^{-1} R_t I_t \\ &-\left(\left(W^{-1} C_t^{-1} + \mathcal{L}_c W\right) L_t^{-1} + W^{-1} C_t^{-1} \mathcal{B} L^{-1} \mathcal{B}^T\right) V \\ &- W^{-1} C_t^{-1} \mathcal{B} L^{-1} R I - G_a u, \end{aligned} \tag{22}$$

and $G_d, G_a \in \mathbb{R}^{n \times n}$ are

$$\begin{aligned} G_d &= (W^{-1} C_t^{-1} + \mathcal{D}_c W) L_t^{-1}, \\ G_a &= \mathcal{A}_c W L_t^{-1}. \end{aligned} \tag{23}$$

Here, $\mathcal{D}_c$ and $\mathcal{A}_c$ are the degree matrix and the adjacency matrix of the communication graph, respectively, i.e. $\mathcal{L}_c = \mathcal{D}_c - \mathcal{A}_c$. We assume that the entries of $b$ and $G_d$ have known bounds for all $i \in \mathcal{V}$:

$$\begin{aligned} |b_i| &\leq b_{\max_i} \\ G_{\min_i} &\leq G_{d_i} \leq G_{\max_i}, \end{aligned} \tag{24}$$

with $b_{\max_i}$, $G_{\min_i}$ and $G_{\max_i}$ being positive constants. According to the theory underlying the so-called Suboptimal

---

[‡] For the sake of simplicity, the order $r$ of the sliding manifold is omitted in the remainder of this paper.

[§] The relative degree is the minimum order $\rho$ of the time derivative $\sigma_i^{(\rho)}, i \in \mathcal{V}$, of the sliding variable associated with the $i$-th node in which the control $u_i, i \in \mathcal{V}$ explicitly appears.

SOSM (SSOSM) control algorithm [30], the $i$-th SOSM control law, that can be used to steer $\xi_{1_i}$ and $\xi_{2_i}$, to zero in a finite time, even in presence of uncertainties, is given by

$$u_i = -\mu_i U_{\max_i} \, \text{sgn} \left( \xi_{1_i} - \tfrac{1}{2} \xi_{1_i}^{\max} \right), \tag{25}$$

with

$$U_{\max_i} > \max \left( \frac{b_{\max_i}}{\mu_i^* G_{\min_i}}; \frac{4 b_{\max_i}}{3 G_{\min_i} - \mu_i^* G_{\max_i}} \right), \tag{26}$$

$$\mu_i^* \in (0,1] \cap \left( 0, \frac{3 G_{\min_i}}{G_{\max_i}} \right), \tag{27}$$

$\mu_i$ switching between $\mu_i^*$ and 1, according to [30, Algorithm 1]. The extremal value $\xi_{1_i}^{\max}$ in (25) can be detected by implementing for instance a peak detector as in [31]. Note that only the value of $\xi_{1_i}$, i.e., $w_i(V_i - V_i^\star) - \theta_i$, is required to generate the control signal $u_i$.

*Remark 5:* (**Switching frequency**) The discontinuous control signal (25) can be directly used in practice to open and close the switch of the Buck converter. As a result, the Insulated Gate Bipolar Transistors (IGBTs) switching frequency cannot be a-priori fixed and the power losses could be high. Usually, in order to achieve a constant IGBTs switching frequency, Buck converters are controlled by implementing the so-called Pulse Width Modulation (PWM) technique. To do this, a continuous control signal, that represents the so-called duty cycle of the Buck converter, is required.

### B. Third Order SM control: duty cycle

To ensure a continuous control input (duty cycle), we adopt the procedure suggested in [30] and first integrate the (discontinuous) control signal generated by a sliding mode controller, yielding for system (3) augmented with (14)

$$\begin{aligned} L_t \dot{I}_t &= -R_t I_t - V + u \\ C_t \dot{V} &= I_t + \mathcal{B} I - I_L \\ L \dot{I} &= -\mathcal{B}^T V - R I \\ \dot{\theta} &= -\mathcal{L}_c W I_t \\ \dot{u} &= v, \end{aligned} \tag{28}$$

where $v$ is the new (discontinuous) control input. Note that the input signal to the converter, $u(t) = \int_0^t v(\tau) d\tau$, is continuous, so that $u_i$ can be used as duty cycle for the switch of the $i$-th Buck converter. A consequence is that the system relative degree (with respect to the new control input $v$) is now equal to three, so that we need to rely on a third order sliding mode (3SM) control strategy to reach the sliding manifold $\{(I_t, V, I, \theta): \sigma = \dot{\sigma} = \ddot{\sigma} = \mathbf{0}\}$ in a finite time. To do so, we define the auxiliary variables $\xi_1 = \sigma$, $\xi_2 = \dot{\sigma}$ and $\xi_3 = \ddot{\sigma}$, and build the auxiliary system as follows

$$\begin{aligned} \dot{\xi}_1 &= \xi_2 \\ \dot{\xi}_2 &= \xi_3 \\ \dot{\xi}_3 &= \dot{b}(I_t, V, I, u) + G_d v \\ \dot{u} &= v, \end{aligned} \tag{29}$$

with $b$ as in (22), and $G_d$, $G_a$ as in (23). Then, we assume that the entries of $\dot{b}$ can be bounded as

$$|\dot{b}_i(\cdot)| \leq \beta_{\max_i} \quad \forall i \in \mathcal{V}, \tag{30}$$

where $\beta_{\max_i}$ is a known positive constant.

*Remark 6:* (**Uncertainty of** $b, \dot{b}$ **and** $G_d$) The mappings $b, \dot{b}$ and matrix $G_d$ are uncertain due to the presence of the unmeasurable current demand $I_L$ and possible network parameter uncertainties. However, relying on Assumption 1 and observing that $b$ and $\dot{b}$ depend on the electric signals related to the finite power of the microgrid, $b, \dot{b}$ and $G_d$ are in practice bounded. Generally, the bounds of the unknown quantities can be determined by data analysis and engineering understanding.

Now, the 3SM control law proposed in [32] can be used to steer $\xi_{1_i}, \xi_{2_i}$ and $\xi_{3_i}, i \in \mathcal{V}$, to zero in a finite time. It is given by

$$v_i = -\alpha_i \begin{cases} v_{1_i} = \text{sgn}(\ddot{\sigma}_i) & \boldsymbol{\sigma}_i \in \mathcal{M}_{1_i}/\mathcal{M}_{0_i} \\ v_{2_i} = \text{sgn}\left( \dot{\sigma}_i + \frac{\ddot{\sigma}_i^2 v_{1_i}}{2\alpha_{r_i}} \right) & \boldsymbol{\sigma}_i \in \mathcal{M}_{2_i}/\mathcal{M}_{1_i} \\ v_{3_i} = \text{sgn}(s_i(\boldsymbol{\sigma}_i)) & \text{otherwise,} \end{cases} \tag{31}$$

where $\boldsymbol{\sigma}_i = [\sigma_i, \dot{\sigma}_i, \ddot{\sigma}_i]^T$ and

$$s_i(\boldsymbol{\sigma}_i) = \sigma_i + \frac{\ddot{\sigma}_i^3}{3\alpha_{r_i}^2} + v_{2_i} \left[ \frac{1}{\sqrt{\alpha_{r_i}}} \left( v_{2_i} \dot{\sigma}_i + \frac{\ddot{\sigma}_i^2}{2\alpha_{r_i}} \right)^{\frac{3}{2}} + \frac{\dot{\sigma}_i \ddot{\sigma}_i}{\alpha_{r_i}} \right],$$

with

$$\alpha_{r_i} = \alpha_i G_{\min_i} - \beta_{\max_i} > 0. \tag{32}$$

Then, given the bounds $G_{\min_i}$ and $\beta_{\max_i}$, the control amplitude $\alpha_i$ is chosen such that $\alpha_{r_i}$ is positive. The manifolds $\mathcal{M}_{1_i}, \mathcal{M}_{2_i}, \mathcal{M}_{3_i}$ in (31) are defined as

$$\mathcal{M}_{0_i} = \{ \boldsymbol{\sigma}_i \in \mathbb{R}^3 : \sigma_i = \dot{\sigma}_i = \ddot{\sigma}_i = 0 \}$$

$$\mathcal{M}_{1_i} = \left\{ \boldsymbol{\sigma}_i \in \mathbb{R}^3 : \sigma_i - \frac{\ddot{\sigma}_i^3}{6\alpha_{r_i}^2} = 0, \dot{\sigma}_i + \frac{\ddot{\sigma}_i |\ddot{\sigma}_i|}{2\alpha_{r_i}} = 0 \right\}$$

$$\mathcal{M}_{2_i} = \{ \boldsymbol{\sigma}_i \in \mathbb{R}^3 : s_i(\boldsymbol{\sigma}_i) = 0 \}.$$

From (31), one can observe that the controller of DGU $i$ requires not only $\sigma_i$, but also $\dot{\sigma}_i$ and $\ddot{\sigma}_i$. Yet, according to Assumption 1, only $I_{t_i}$ and $V_i$ are measurable at the $i$-th DGU. Then, one can rely on Levant's second-order differentiator [33] to retrieve $\dot{\sigma}_i$ and $\ddot{\sigma}_i$ in a finite time. Consequently, for system (29), the estimators are given by

$$\begin{aligned} \dot{\hat{\xi}}_{1_i} &= -\lambda_{0_i} \left| \hat{\xi}_{1_i} - \xi_{1_i} \right|^{\frac{2}{3}} \text{sgn}\left( \hat{\xi}_{1_i} - \xi_{1_i} \right) + \hat{\xi}_{2_i} \\ \dot{\hat{\xi}}_{2_i} &= -\lambda_{1_i} \left| \hat{\xi}_{2_i} - \dot{\hat{\xi}}_{1_i} \right|^{\frac{1}{2}} \text{sgn}\left( \hat{\xi}_{2_i} - \dot{\hat{\xi}}_{1_i} \right) + \hat{\xi}_{3_i} \\ \dot{\hat{\xi}}_{3_i} &= -\lambda_{2_i} \, \text{sgn}\left( \hat{\xi}_{3_i} - \dot{\hat{\xi}}_{2_i} \right), \end{aligned} \tag{33}$$

where $\hat{\xi}_{1_i} = \hat{\sigma}_i$, $\hat{\xi}_{2_i} = \dot{\hat{\sigma}}_i$ and $\hat{\xi}_{3_i} = \ddot{\hat{\sigma}}_i$ are the estimated values of $\xi_{1_i} = \sigma_i, \xi_{2_i} = \dot{\sigma}_i$ and $\xi_{3_i} = \ddot{\sigma}_i$, respectively. The estimates obtained via (33) can be used in (31), replacing the original variables. The other parameters are $\lambda_{0_i} = 3\Lambda_i^{1/3}$, $\lambda_{1_i} = 1.5\Lambda_i^{1/2}$, $\lambda_{2_i} = 1.1\Lambda_i$, $\Lambda_i > 0$, as suggested in [33].

*Remark 7:* (**Scalability and distributed control**) Since the selected sliding function (20) is designed by using the additional state $\theta$ in (14), the overall control scheme is indeed distributed, and only information on generated currents $I_t$ needs to be shared. More precisely, the controller of the $i$-th DGU needs information only from the DGUs that communicate with it. Note that the design of the local controller for each DGU is not based on the knowledge of the whole microgrid, so that the complexity of the control synthesis does not depend on the microgrid size.

*Remark 8:* (**Alternative SM controllers**) In this work we rely on the SOSM control algorithm proposed in [30] and the 3SM control law proposed in [32]. However, the results in this paper are obtained independent of the particular choice of sliding mode controller.

## VI. Stability Analysis

In this section we first show that the states of the controlled microgrid are constrained, after a finite time, to the manifold $\sigma = \mathbf{0}$, where Objective 2 is achieved. Thereafter, we prove that the solutions to the system, once the sliding manifold is attained, converge exponentially to a constant point, achieving additionally Objective 1.

### A. Equivalent reduced order system

As a first step, we study the convergence to the sliding manifold when the SSOSM or the 3SM control law is applied to the system.

*Lemma 4:* (**Convergence to the sliding manifold: SSOSM**) Let Assumption 1 hold. The solutions to system (3) augmented with (14), controlled via the SSOSM control law (25), converge in a finite time $T_r$, to the sliding manifold $\{(I_t, V, I, \theta) : \sigma = \dot{\sigma} = \mathbf{0}\}$, with $\sigma$ given by (20).

*Proof:* Following [30], the application of (25) to each converter guarantees that $\sigma = \dot{\sigma} = \mathbf{0}$, for all $t \geq T_r$. ∎

*Lemma 5:* (**Convergence to the sliding manifold: 3SM**) Let Assumption 1 hold. The solutions to system (3) augmented with (14), controlled via 3SM control algorithm (29)-(33), converge in a finite time $T_r$, to the sliding manifold $\{(I_t, V, I, \theta) : \sigma = \dot{\sigma} = \ddot{\sigma} = \mathbf{0}\}$, with $\sigma$ given by (20).

*Proof:* By implementing the Levant's differentiator (33) in each node, the values of $\xi_1, \xi_2, \xi_3$, are estimated in a finite time $T_{Ld} \geq 0$ [33]. Then, following [32], the application of (31) to each converter guarantees that $\sigma = \dot{\sigma} = \ddot{\sigma} = \mathbf{0}$, for all $t \geq T_r \geq T_{Ld}$. ∎

As we will show in the proof of Theorem 2 in the next subsection, converging to the sliding manifold where $\sigma = \mathbf{0}$, is sufficient to conclude that Objective 2 (voltage balancing) is achieved. We postpone the analysis, in order to show additionally convergence to a *constant* voltage. For the analysis of the system, when the solutions are constrained to the sliding manifold, it is convenient to exploit the so-called system order reduction property, typical of sliding mode control methodology. Indeed, when the state of system (3) augmented with (14) is constrained to the sliding manifold $\{(I_t, V, I, \theta) : \sigma = \dot{\sigma} = \mathbf{0}\}$, with $\sigma$ given by (20), the controlled system is described by $3n+m$ differential equations and

$2n$ algebraic equations. Then, it is possible to obtain $2n$ state variables depending on the other $n+m$ ones. The resulting system of order $n+m$ represents the reduced order system equivalent to the system controlled with a discontinuous law, with the initial condition $(I_t(T_r), V(T_r), I(T_r), \theta(T_r))$, when $\sigma = \dot{\sigma} = \mathbf{0}$.

*Lemma 6:* (**Equivalent reduced order system**) For all $t \geq T_r$, the dynamics of the controlled system (3) augmented with (14) are given by the following equivalent system of reduced order

$$
\begin{aligned}
C_t \dot{V} &= \left(\mathbb{I}_n - (\mathbb{I}_n + C_t W \mathcal{L}_c W)^{-1}\right) \mathcal{B} I \\
&\quad - \left(\mathbb{I}_n - (\mathbb{I}_n + C_t W \mathcal{L}_c W)^{-1}\right) I_L \\
L\dot{I} &= -\mathcal{B}^T V - RI,
\end{aligned}
\tag{34}
$$

together with the following algebraic relations

$$
\theta = W^{-1} (V - V^\star)
\tag{35}
$$

$$
I_t = (\mathbb{I}_n + C_t W \mathcal{L}_c W)^{-1} (-\mathcal{B}I + I_L).
\tag{36}
$$

*Proof:* Given the sliding function (20), by virtue of Lemma 4 and Lemma 5, the state of system (3) augmented with (14) is constrained to the manifold $\{(I_t, V, I, \theta) : \sigma = \dot{\sigma} = \mathbf{0}\}$, where $\theta = W^{-1} (V - V^\star)$ and $\dot{V} = W\dot{\theta}$. From the latter, one can straightforwardly obtain (36). After substituting expression (36) for $I_t$ in (3), the dynamics of the voltage $V$ become as in (34). ∎

### B. Exponential convergence and objectives attainment

In the pervious subsection, we established that after a finite time $T_r$, the dynamics of the controlled microgrid are described by the equivalent system (34). In this subsection we study the convergence properties of this equivalent system. To do so, we rely on the concept of semistability [34], of which we recall the definition for convenience.

*Definition 4:* (**Semistability**) Consider the autonomous system

$$
\dot{x}(t) = Ax(t),
\tag{37}
$$

where $t \geq 0$, $x \in \mathbb{R}^n$ and $A \in \mathbb{R}^{n \times n}$. System (37) is semistable if $\lim_{t \to \infty} x(t)$ exists for all initial conditions $x(0)$.

Furthermore, the following lemma turns out to be useful in the upcoming analysis:

*Lemma 7:* Given a positive definite matrix $P \in \mathbb{R}^{n \times n}$ and a positive semidefinite matrix $Q \in \mathbb{R}^{n \times n}$, then

$$
P - \left(P^{-1} + Q\right)^{-1} \succeq 0.
\tag{38}
$$

*Proof:* Let $\tilde{Q} = P^{\frac{1}{2}} Q P^{\frac{1}{2}}$. Clearly, $\tilde{Q} \succeq 0$, and $\mathbb{I}_n + \tilde{Q} \succ 0$. Then,

$$
P - \left(P^{-1} + Q\right)^{-1} = P^{\frac{1}{2}} \left[\mathbb{I}_n - \left(\mathbb{I}_n + \tilde{Q}\right)^{-1}\right] P^{\frac{1}{2}}
\tag{39}
$$

is a positive semidefinite matrix if and only if $\mathbb{I}_n - (\mathbb{I}_n + \tilde{Q})^{-1} = \tilde{Q}(\mathbb{I}_n + \tilde{Q})^{-1} \succeq 0$. Observing that $(\mathbb{I}_n + \tilde{Q})^{-1} \succ 0$, it yields

$$
\tilde{Q}(\mathbb{I}_n + \tilde{Q})^{-1} \backsim (\mathbb{I}_n + \tilde{Q})^{-\frac{1}{2}} \tilde{Q}(\mathbb{I}_n + \tilde{Q})^{-\frac{1}{2}} \succeq 0,
\tag{40}
$$

which completes the proof. ∎

Next, we show that the line currents $I$ converge to a constant value.

*Lemma 8:* (**Convergence of** $I$) Let Assumptions 1 and 2 hold. Given the equivalent reduced order system (34), $\lim_{t \to \infty} I(t)$ exists for all initial conditions $I(T_r)$.

*Proof:* Let $\tilde{V} = V - \overline{V}$ and $\tilde{I} = I - \overline{I}$ be the error given by the difference between the state of system (34) and a steady state value. Then, the dynamics of the corresponding error system are given by

$$C_t \dot{\tilde{V}} = \left( \mathbb{1}_n - (\mathbb{1}_n + C_t W \mathcal{L}_c W)^{-1} \right) \mathcal{B} \tilde{I}$$
$$L \dot{\tilde{I}} = -\mathcal{B}^T \tilde{V} - R\tilde{I}, \qquad (41)$$

From (41), we obtain

$$\dot{\tilde{V}} = C_t^{-1} \left( \mathbb{1}_n - (\mathbb{1}_n + C_t W \mathcal{L}_c W)^{-1} \right) \mathcal{B} \tilde{I}$$
$$= \left( C_t^{-1} - (C_t + C_t W \mathcal{L}_c W C_t)^{-1} \right) \mathcal{B} \tilde{I}, \qquad (42)$$

and

$$L \ddot{\tilde{I}} + R\dot{\tilde{I}} + \mathcal{B}^T \dot{\tilde{V}} = \mathbf{0}. \qquad (43)$$

Substituting expression (42) for $\dot{\tilde{V}}$ in (43) leads to

$$L\ddot{\tilde{I}} + R\dot{\tilde{I}} + \underbrace{\mathcal{B}^T \left( C_t^{-1} - (C_t + C_t W \mathcal{L}_c W C_t)^{-1} \right) \mathcal{B}}_{K} \tilde{I} = \mathbf{0}. \qquad (44)$$

Since, by virtue of Lemma 7 (with $P = C_t^{-1}, Q = C_t W \mathcal{L}_c W C_t$), $C_t^{-1} - (C_t + C_t W \mathcal{L}_c W C_t)^{-1} \succeq 0$, then we also have $K = \mathcal{B}^T (C_t^{-1} - (C_t + C_t W \mathcal{L}_c W C_t)^{-1}) \mathcal{B} \succeq 0$. According to [34, Corollary 2], system (44) is semistable (see Definition 4) if and only if

$$\text{rank} \begin{bmatrix} R \\ R(L^{-1}K) \\ R(L^{-1}K)^2 \\ \vdots \\ R(L^{-1}K)^{m-1} \end{bmatrix} = m. \qquad (45)$$

Since $R$ is a positive definite $m \times m$ diagonal matrix it can be readily confirmed that condition (45) holds, such that system (44) is indeed semistable. Since $\overline{I}$ is a constant vector, it immediately follows that $\lim_{t \to \infty} I(t)$ exists. ∎

Lemma 8 established that $\lim_{t \to \infty} I(t)$ exists for all initial conditions $I(T_r)$. This result can now be exploited to show that also the voltages converge to constant values.

*Lemma 9:* (**Convergence of** $V$) Let Assumptions 1–3 hold. Given the equivalent reduced order system (34), $\lim_{t \to \infty} V(t)$ exists for all initial conditions $V(T_r)$.

*Proof:* Exploiting the convergence of $I$ to a constant vector (see Lemma 8), from (43) we have

$$\lim_{t \to \infty} \mathcal{B}^T \dot{V}(t) = \mathbf{0}, \qquad (46)$$

implying that

$$\lim_{t \to \infty} \dot{V}(t) = \mathbb{1}_n \kappa, \qquad (47)$$

with $\kappa \in \mathbb{R}$. By virtue of Lemma 3 and Lemma 4 or Lemma 5, for all $t \geq T_r$, we also have

$$\mathbb{1}_n^T W^{-1} V = \mathbb{1}_n^T (\theta + W^{-1} V^\star) = \mathbb{1}_n^T \theta(0) + \mathbb{1}_n^T W^{-1} V^\star. \qquad (48)$$

Taking the derivative with respect to time on both sides of (48), it follows that $\mathbb{1}_n^T W^{-1} \dot{V}(t) = 0$ for all $t \geq T_r$. Exploiting (47), we obtain

$$\lim_{t \to \infty} \mathbb{1}_n^T W^{-1} \dot{V}(t) = \mathbb{1}_n^T W^{-1} \lim_{t \to \infty} \dot{V}(t)$$
$$= \mathbb{1}_n^T W^{-1} \mathbb{1}_n \kappa \qquad (49)$$
$$= 0,$$

which implies $\kappa = 0$ and consequently that $\lim_{t \to \infty} V(t)$ exists for all initial conditions $V(T_r)$. ∎

We are now ready to establish the first main result of this paper.

*Theorem 1:* (**Achieving current sharing**) Let Assumptions 1–3 hold. Consider system (3), (14), controlled with the proposed distributed SSOSM (Subsection V-A) or 3SM (Subsection V-B) control scheme. Then, the controlled currents $I_t(t)$ converge, after a finite time, exponentially to $W^{-1} \mathbb{1}_n \mathbb{1}_n^T I_L / (\mathbb{1}_n^T W^{-1} \mathbb{1}_n)$, achieving proportional current sharing.

*Proof:* According to Lemma 6, for all $t \geq T_r$, the dynamics of the controlled system (3), (14) are given by the autonomous system (34) together with the algebraic equations (35) and (36). Bearing in mind the results proved in Lemma 8 and Lemma 9, the dynamics of the line current $I$ and the voltage $V$ are semistable. From the algebraic equations (35) and (36), it follows that $\lim_{t \to \infty} \theta(t)$ and $\lim_{t \to \infty} I_t(t)$ exist as well. Since (34) is linear and $\ker(\mathcal{L}_c) = \text{im}(\mathbb{1}_n)$, (14) implies that the vector $I_t(t)$, with initial condition $I_t(T_r)$, converges exponentially to a constant vector, achieving proportional current sharing. ∎

We now proceed with establishing the second main result of this paper.

*Theorem 2:* (**Achieving voltage balancing**) Let Assumptions 1–3 hold. Consider system (3), (14), controlled with the proposed distributed SSOSM (Subsection V-A) or 3SM (Subsection V-B) control scheme. Then, given a desired references vector $V^\star$, the voltages $V(t)$ satisfy $\mathbb{1}_n^T W^{-1} V(t) = \mathbb{1}_n^T W^{-1} V^\star$ for all $t \geq T_r$, with $T_r$ a finite time.

*Proof:* Following Lemma 4 or Lemma 5, for all $t \geq T_r$, the equality $W^{-1} V(t) = W^{-1} V^\star + \theta(t)$ holds. Pre-multiplying both sides by $\mathbb{1}_n^T$ yields $\mathbb{1}_n^T W^{-1} V(t) = \mathbb{1}_n^T W^{-1} V^\star + \mathbb{1}_n^T \theta(t)$. Due to Assumption 3 and by virtue of Lemma 3, one has that $\mathbb{1}_n^T \theta(t) = \mathbb{1}_n^T \theta(0) = 0$. Then, one can conclude that voltage balancing is achieved for all $t \geq T_r$. ∎

*Remark 9:* (**Robustness to failed communication**) The proposed control scheme is distributed and as such requires a communication network to share information on the generated currents. However, note that the integrators $\theta$ in (14) are not needed to regulate the voltages in the microgrid to their desired values, but are only required to achieve current sharing and voltage balancing. In fact, by omitting the variable $\theta$ in the analysis, the controlled microgrid converges, in a finite time, to the manifold $\sigma = \mathbf{0}$, where $V = V^\star$, as shown in [14]. Moreover, considering constant value of $\theta_i$ (e.g. after the plug-out of the DGU $i$), the controlled DGU $i$ converges, in a finite time, to the manifold $\sigma_i = 0$, where $V_i = V_i^\star + w_i \overline{\theta}_i$.

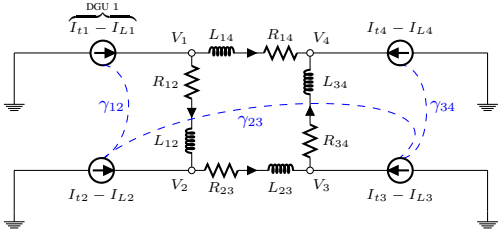

Fig. 2. Scheme of the considered (Kron reduced) microgrid with 4 power converters. The dashed lines represent the communication network.

TABLE II
MICROGRID PARAMETERS AND CURRENT DEMAND

| DGU | | 1 | 2 | 3 | 4 |
|---|---|---|---|---|---|
| $R_{ti}$ | (Ω) | 0.2 | 0.3 | 0.5 | 0.1 |
| $L_{ti}$ | (mH) | 1.8 | 2.0 | 3.0 | 2.2 |
| $C_{ti}$ | (mF) | 2.2 | 1.9 | 2.5 | 1.7 |
| $w_i$ | (–) | $0.4^{-1}$ | $0.2^{-1}$ | $0.15^{-1}$ | $0.25^{-1}$ |
| $V_i^\star$ | (V) | 380.0 | 380.0 | 380.0 | 380.0 |
| $I_{Li}(0)$ | (A) | 30.0 | 15.0 | 30.0 | 26.0 |
| $\Delta I_{Li}$ | (A) | 10.0 | 7.0 | −10.0 | 5.0 |

TABLE III
LINE PARAMETERS

| Line | | {1,2} | {2,3} | {3,4} | {1,4} |
|---|---|---|---|---|---|
| $R_{ij}$ | (mΩ) | 70 | 50 | 80 | 60 |
| $L_{ij}$ | (μH) | 2.1 | 2.3 | 2.0 | 1.8 |

*Remark 10:* (**Perturbations in the controller states**) In case Assumption 3 is violated, we have $\mathbb{1}_n^T \theta(t) = \mathbb{1}_n^T \theta(0)$, and consequently $\mathbb{1}_n^T W^{-1} V(t) = \mathbb{1}_n^T W^{-1} V^\star + \mathbb{1}_n^T \theta(0)$ on the sliding manifold, implying that the weighted average voltage of the microgrid is shifted by $\mathbb{1}_n^T \theta(0)$. However, the presented stability analysis is still valid such that the stability of the whole microgrid and the achievement of proportional current sharing is still guaranteed.

## VII. SIMULATION RESULTS

In this section, the proposed manifold-based consensus algorithm is assessed in simulation by implementing the third order sliding mode control strategy discussed in Subsection V-B. We consider a microgrid composed of 4 DGUs interconnected as shown in Figure 2, where also the communication network is depicted. The parameters of each DGU, including the current demand, and the line parameters are reported in Tables II and III, respectively. The weights associated with the edges of the communication graph are $\gamma_{12} = \gamma_{23} = \gamma_{34} = 10$. For all the DGUs the controller parameter $\alpha_i$ in (31) is set to $2.4 \times 10^3$. In order to investigate the performance of the proposed control approach within a low voltage DC microgrid, four different scenarios are implemented (see Fig. 3).

### A. Scenario 1: proportional current sharing

The system is initially at the steady state. Then, consider a current demand variation $\Delta I_{Li}$ at the time instant $t = 1$ s (see Table II). The PCC voltages and the generated currents are illustrated in Figure 4. One can appreciate that the weighted average of the PCC voltages (denoted by $V_{av}$) is always equal to the weighted average of the corresponding references (see Objective 2), and the current generated by each DGU converges to the desired value, achieving proportional current sharing (see Objective 1). Moreover, in Figure 5 the currents shared among the DGUs are reported together with the control signals generated by the 3SM algorithm (31). Note that the 3SM controllers, which require only local measurements of $V_i$ and information on $I_t$ from neighbours over the communication network, generate control signals that are equal to the optimal feedforward input (7), without exact knowledge on the network parameters and the current demand $I_L$.

### B. Scenario 2: opening of a distribution line

In the second scenario, we investigate the performance of the proposed controllers when a distribution line is opened (e.g. due to an electric fault). The system is initially at the steady state, and at the time instant $t = 0.4$ s, the distribution line interconnecting the DGUs 1 and 4 is opened. Then, consider a current demand variation as in Scenario 1. The PCC voltages and the generated currents are illustrated in Figure 6. One can appreciate that the weighted average of the PCC voltages (denoted by $V_{av}$) is always equal to the weighted average of the corresponding references (see Objective 2), and the current generated by each DGU converges to the desired value, achieving proportional current sharing (see Objective 1). Moreover, in Figure 7 the currents shared among the DGUs are reported together with the control signals generated by the 3SM control algorithm (31).

### C. Scenario 3: plug-out and plug-in of a DGU

In the third scenario, we investigate the Plug-and-Play (PnP) capabilities of the proposed controllers. For the sake of clarity, in this scenario and the next one we consider equal current sharing among the DGUs. The system is initially at the steady state, and at the time instant $t = 0.4$ s, the DGU 4 is disconnected from the considered DC network (in this configuration the impedance of the line interconnecting DGU 1 and DGU 3 is equal to the sum of the line impedances $Z_{14}$ and $Z_{34}$). After a current demand variation as in Scenario 1, at the time instant $t = 1.4$ s, the DGU 4 is reconnected to the DC network. The PCC voltages and the generated currents are illustrated in Figure 8. One can appreciate that the arithmetic average of the PCC voltages (denoted by $V_{av}$) is equal to the arithmetic average of the corresponding references, even after disconnecting the DGU 4. Moreover, when the DGU 4 operates isolated from the considered DC network, equal current sharing is achieved only among the DGUs 1, 2 and 3, while the DGU 4 supplies its local load. However, when the DGU 4 is reconnected to the DC network, current sharing among all the DGUs is again reestablished. Moreover, in Figure 9 the currents shared among the DGUs are reported together with the control signals generated by the 3SM control algorithm (31). Note that, when the DGU 4 is isolated from the network, the comparison between $u_4$ and the corresponding optimal feedforward input loses its meaning.

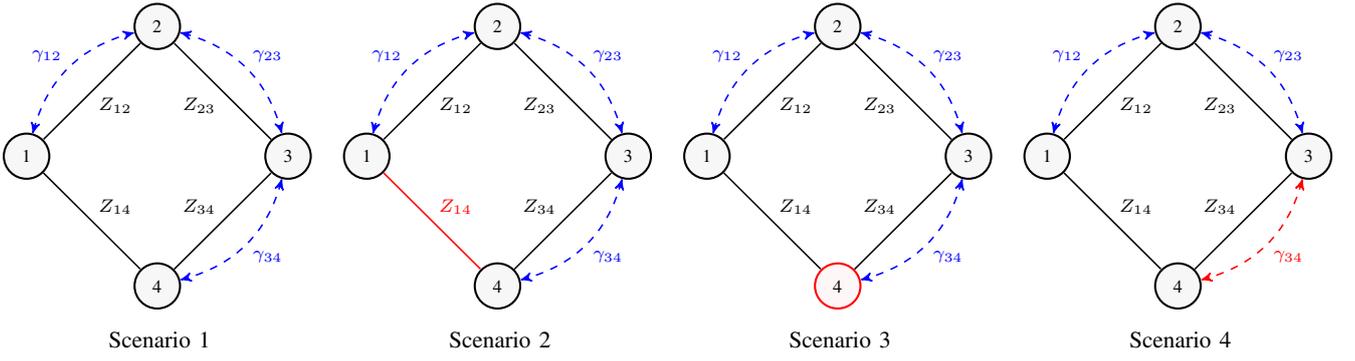

Fig. 3. From the left: the configurations of the considered microgrid implemented in Scenario 1, Scenario 2, Scenario 3 and Scenario 4, respectively. $Z_{ij}$ denotes the resistive-inductive impedance of the distribution line interconnecting DGU $i$ with DGU $j$. Components that are failing/removed during the simulation are colored red.

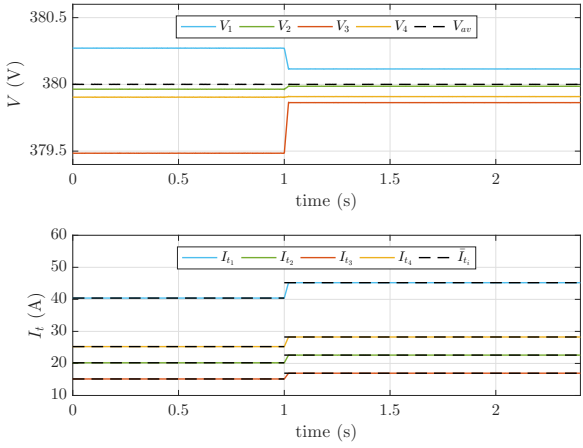

Fig. 4. Scenario 1. From the top: voltage at the PCC of each DGU together with its weighted average value (dashed line); generated currents together with the corresponding values (dashed lines) that allow to achieve proportional current sharing.

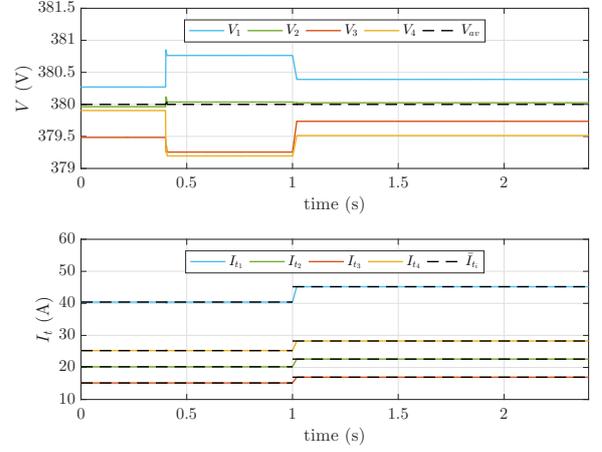

Fig. 6. Scenario 2. From the top: voltage at the PCC of each DGU together with its weighted average value (dashed line); generated currents together with the corresponding values (dashed lines) that allow to achieve proportional current sharing.

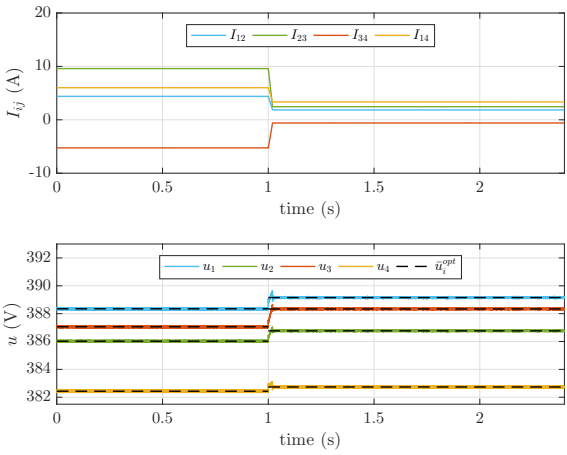

Fig. 5. Scenario 1. From the top: currents shared among the DGUs through the lines; control inputs $u_i(t) = \int_0^t v_i(\tau)d\tau$, $v_i$ as in (31), together with the optimal feedforward inputs (7) indicated by the dashed lines.

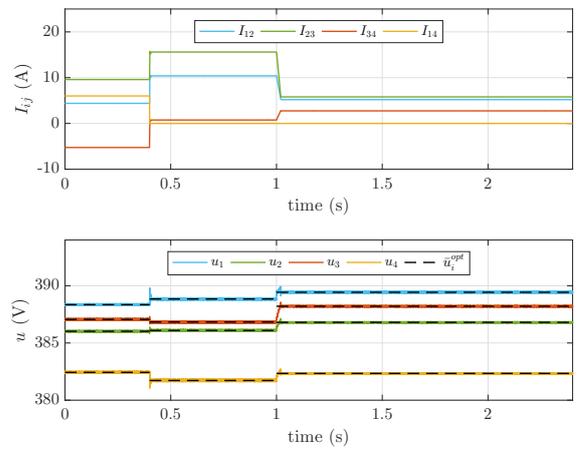

Fig. 7. Scenario 2. From the top: currents shared among the DGUs through the lines; control inputs $u_i(t) = \int_0^t v_i(\tau)d\tau$, $v_i$ as in (31), together with the optimal feedforward inputs (7) indicated by the dashed lines.

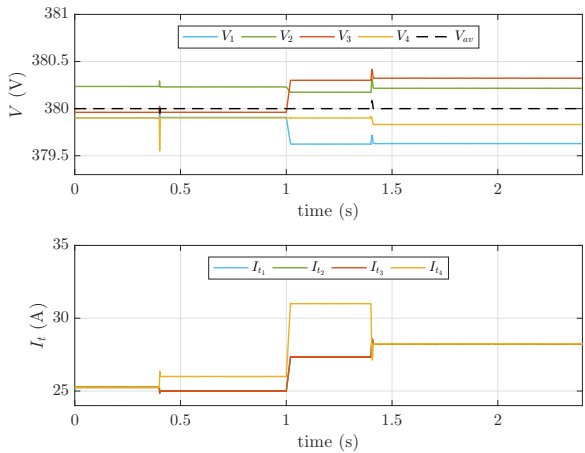

Fig. 8. Scenario 3. From the top: voltage at the PCC of each DGU together with its average value (dashed line); generated currents in case of equal current sharing, which is achieved by DGUs 1, 2, 3 for all the simulation time interval, and by DGU 4 only when it is connected to the microgrid.

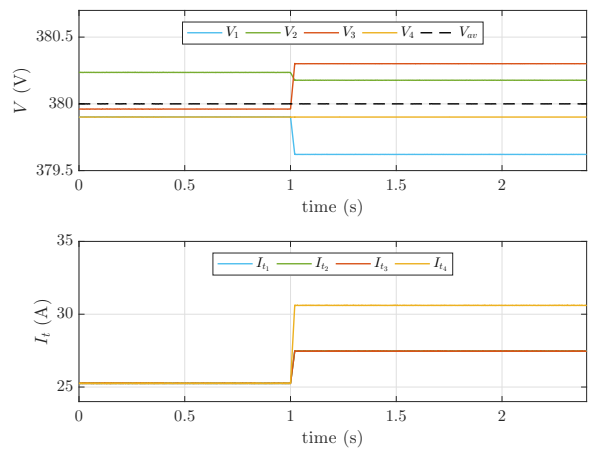

Fig. 10. Scenario 4. From the top: voltage at the PCC of each DGU together with its average value (dashed line); generated currents in case of equal current sharing, which is achieved by DGUs 1, 2, 3 for all the simulation time interval, and by DGU 4 only before the failing of the communication link.

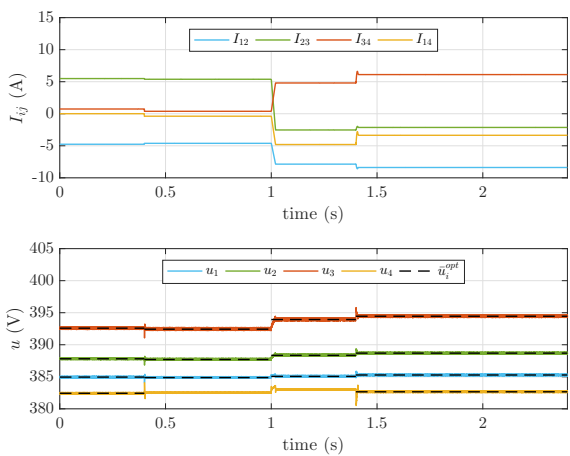

Fig. 9. Scenario 3. From the top: currents shared among the DGUs through the lines; control inputs $u_i(t) = \int_0^t v_i(\tau)d\tau$, $v_i$ as in (31), together with the optimal feedforward inputs (7) indicated by the dashed lines.

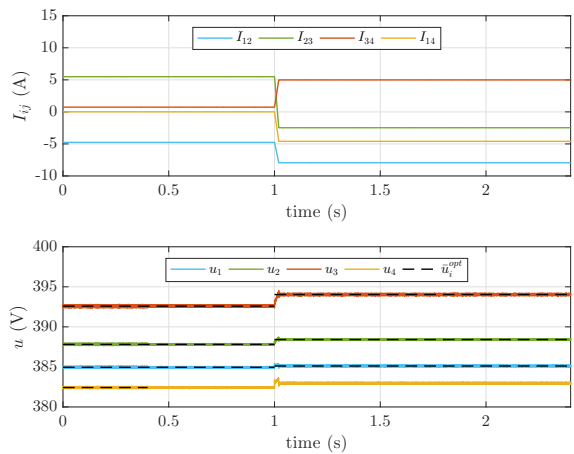

Fig. 11. Scenario 4. From the top: currents shared among the DGUs through the lines; control inputs $u_i(t) = \int_0^t v_i(\tau)d\tau$, $v_i$ as in (31), together with the optimal feedforward inputs (7) indicated by the dashed lines.

### D. Scenario 4: failing of a communication link

In the last scenario, we investigate the robustness of the proposed controllers to failed communication. The system is initially at the steady state, and at the time instant $t = 0.4\,\text{s}$, the communication between DGU 3 and DGU 4 is interrupted. We observe that as long as the demand does not change, current sharing among all the DGUs in mainteined. The PCC voltages and the generated currents are illustrated in Figure 10. One can note that after a current demand variation (see Table II), equal current sharing is achieved only among the DGUs 1, 2 and 3, while the DGU 4 generates a current such that the voltage at node 4 is kept constant. One can appreciate that the arithmetic average of the PCC voltages (denoted by $V_{av}$) is equal to the arithmetic average of the corresponding references, even after interrupting the communication between DGU 3 and DGU 4. Moreover, in Figure 11 the currents shared among the DGUs are reported together with the control signals generated by the 3SM control algorithm (31). Note that, when the communication with DGU 4 fails, the comparison between $u_4$ and the corresponding optimal feedforward input loses its meaning.

Even if IEEE Standards or guidlines for DC power distribution networks do not exist yet (to the best of our knowledge), it is usually required in practical cases that the voltage deviations are within the 5 % of the desired value (see for instance [24]). In all the previous scenarios, the voltage at the PCC of each DGU is within the range $380 \pm 1$ V, implying that the voltage deviations are less than the 0.3 % of the nominal value $V^\star = 380$ V, even during transients and critical conditions.

## VIII. Conclusions

In this paper we design a distributed control algorithm, obtaining proportional current sharing and voltage regulation in DC microgrids. Its convergence properties are analytically investigated. The proposed control scheme exploits a communication network to achieve proportional current sharing using a consensus-like algorithm. Another useful feature of the proposed control scheme is that the weighted average voltage of the microgrid converges to the weighted average of the voltage references, independently of the initial voltage conditions. The latter is achieved by constraining the system to a suitable manifold. To ensure that the desired manifold is reached in a finite time, even in presence of modelling uncertainties, two sliding mode control strategies are proposed, that provide the switching frequencies or the duty cycle of the power converters. An extensive simulation analysis is also provided, considering different and realistic scenarios. The proposed controllers show satisfactory closed-loop performance and Plug-and-Play capabilities, even in presence of uncertainties, topology changes and communication failings. Interesting future research includes the design of distributed controllers aimed at guaranteeing power sharing and studying the stability analysis of the resulting nonlinear system.